%% file: text3.tex
\def\draft{n}
\documentclass[12pt]{amsart}
\usepackage{fullpage,amssymb,epic,eepic}

\theoremstyle{plain}

\newtheorem{theorem}{Theorem}
\newtheorem{proposition}{Proposition}[section]
\newtheorem{lemma}[proposition]{Lemma}
\newtheorem{corollary}[proposition]{Corollary}
\newtheorem{claim}[proposition]{Claim}

\theoremstyle{definition}

\theoremstyle{remark}

\newtheorem{remark}[proposition]{Remark}

\def\printname#1{
	\if\draft y
		\smash{\makebox[0pt]{\hspace{-0.5in}
			\raisebox{8pt}{\tt\tiny #1}}}
	\fi
}

\newlength{\standardunitlength}
\catcode`\@=11
\long\def\@makecaption#1#2{%
    \vskip 10pt
    \setbox\@tempboxa\hbox{
      \small\sf{\bfcaptionfont #1. }\ignorespaces #2}%
    \ifdim \wd\@tempboxa >\captionwidth {%
        \rightskip=\@captionmargin\leftskip=\@captionmargin
        \unhbox\@tempboxa\par}%
      \else
        \hbox to\hsize{\hfil\box\@tempboxa\hfil}%
    \fi}
\font\bfcaptionfont=cmssbx10 scaled \magstephalf
\newdimen\@captionmargin\@captionmargin=2\parindent
\newdimen\captionwidth\captionwidth=\hsize
\catcode`\@=12

\def\lbl#1{\label{#1}\printname{#1}}

\def\eqdef{\overset{\text{def}}{=}}

\def\eqdef{\overset{\text{def}}{=}}

\def\BZ{\mathbb Z}

\def\L{\mathcal L}

\def\s{\sigma}

\def\l{\lambda}

\def\S{\Sigma}

\def\b{\beta}
\def\e{\epsilon}
\def\a{\alpha}
\def\g{\gamma}

\def\la{\langle}
\def\ra{\rangle}

\def\sgn{\text{sgn}'}
\def\n{\nu}
\def\det{\text{det}}
\def\tByt#1#2#3#4{ \begin{pmatrix}
                   #1 & #2 \\
                   #3 & #4 
                   \end{pmatrix} }

\begin{document}


\title[Signatures of links and finite type invariants of
      cyclic branched covers]{Signatures of 
      links and finite type invariants of
      cyclic branched covers}

\author{Stavros Garoufalidis}
\address{Department of Mathematics \\
         Harvard University \\
         1 Oxford Street \\
         Cambridge,  MA 02138, USA}
\email{stavros@math.harvard.edu}
\thanks{The  author was partially supported by an NSF grant 
       DMS-98-00703.\newline
       This and related preprints can also be obtained at
{\tt http://www.math.brown.edu/$\sim$stavrosg } \newline
1991 {\em Mathematics Classification.} Primary 57N10. Secondary 57M25.
}

\dedicatory{Dedicated to Mel Rothenberg.}

\date{This edition: September 1, 1998; First edition: November 10, 1997  }

\begin{abstract}
Recently, 
Mullins calculated  the Casson-Walker invariant of
the 2-fold cyclic branched cover of an oriented link in $S^3$
in terms of its Jones polynomial and its signature, under the assumption
that the 2-fold branched cover is a rational homology 3-sphere.
Using elementary principles, we provide a similar calculation for the general
case. In addition, we calculate the LMO invariant of the $p$-fold
branched cover of twisted knots in $S^3$ in terms of the Kontsevich integral
of the knot. 
\end{abstract}

\maketitle

\tableofcontents


\section{Introduction}
\lbl{sec.intro} 

Given an oriented link $L$ in (oriented) $S^3$, one can associate to it
a family of (oriented) 3-manifolds, namely its $p$-fold {\em 
cyclic branched covers} 
$\S^p_L$, where $p$ is a positive integer. Using these 3-manifolds, one
can associate a family of integer-valued invariants of the link $L$, 
namely its $p$-{\em signatures},
$\s_p(L)$. These signatures, being concordance invariants,
play a key  role in the approach
to link theory via surgery theory.  
  
On the other hand, any numerical invariant of 3-manifolds, evaluated
at the p-fold branched cover, gives numerical invariants of oriented links.
The seminal ideas of mathematical physics, initiated by Witten \cite{Wi}
have recently
produced two axiomatizations (and constructions) of numerical invariants of 
links and 3-manifolds; one under the name of topological quantum field theory
(e.g. \cite{At, RT1, RT2}) and another under the name of finite
type invariants (e.g. \cite{Oh, LMO, BGRT}). Moreover, each of these two 
approaches
 offers a conceptual unification of previously known numerical invariants 
of links and 3-manifolds, such as the {\em Casson invariant} \cite{AM} 
and the {\em Jones 
polynomial} \cite{J}.

It turns out that the Casson invariant $\l$, extended
to all rational homology 3-spheres by Walker \cite{Wa}, and further
extended to all 3-manifolds by Lescop \cite{Le} equals (up to a sign)
to twice  
the degree 1 part of a graph-valued invariant of 3-manifolds \cite{LMO, LMMO}
which turns out to be a 
universal finite type invariant of integral homology 3-spheres, \cite{L}.

Recently D. Mullins \cite{Mu} discovered a beautiful relation between 
the value of 
Casson-Walker invariant of the 2-fold branched cover of a link $L$ in $S^3$
in terms of the link's (2)-signature $\s(L)=\s_2(L)$
and the value of its Jones polynomial $J_L$  
at $-1$, under the assumption that the 2-fold branched cover is a rational
homology 3-sphere. It is a natural question to ask whether this assumption
is really needed.
We can now state our result (where the nullity $\n(L)$ of a link $L$ 
can be defined as the first betti number of $\S^2_L$):

\begin{theorem}
\lbl{thm.1}
For an oriented link $L$ in $S^3$ we have:
\begin{equation}
\lbl{eq.blah}
i^{\s(L)+ \n(L)} \l(\S^2_L) 
= \frac{1}{6} J'_L(-1) + \frac{1}{4}J_L(-1) \s(L)
\end{equation}
\end{theorem}

A few remarks are in order:
\begin{remark}
\lbl{rem.mullins}
In case $\S^2_L$ is a rational homology 3-sphere, the above formula
is Mullin's theorem, as expected. See also Remark \ref{rem.vj}.
\end{remark}

\begin{remark}
For the class of links such that $\S^2_L$ is a rational 
homology 3-sphere (e.g. for all knots), there is a skein theory
relation of the signature, see \cite{Li,Mu}. However, the literature
on signatures 
seems to be avoiding the rest of the links. Our result shows that such
a restriction is not necessary.
\end{remark}

\begin{remark}
Mullin's proof uses the Kauffman bracket definition of the Jones polynomial
\cite{Ka}, and oriented, as well as unoriented, smoothings of the link
that are special to the Jones polynomial. Our proof, simpler and derived from 
first principles,  does not use any of
these special properties of the Jones polynomial. In addition, it provides 
a hint for general relations between link signatures and finite type link and 
3-manifold invariants.
\end{remark}

\begin{corollary}
\lbl{cor.1}
If $L$ is a link with nullity at least 4, then $J_L(-1)=J'_L(-1)=0$.
\end{corollary}

It is natural to ask whether Theorem \ref{thm.1} can be extended to the
case of more general covers (such as $p$-fold branched covers), as well
as the case of more general 3-manifold invariants, such as the LMO invariant
$Z^{LMO}$, or its degree at most $n$ part $Z^{LMO}_n$.
In this direction, we have the following partial result. Let $D_mK$
denote the $m$-fold twisted double of a knot $K$ in $S^3$, see Figure
\ref{double}.

\begin{figure}[htpb]
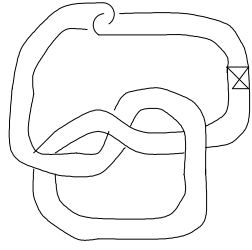

$$ \printname{double}
	\setlength{\unitlength}{0.03\standardunitlength}
	\begin{array}{c}  \hspace{-1.7mm}
        	\raisebox{-8pt}{\input draws/double.tex }
        	\hspace{-1.9mm}
	\end{array}
 $$
\caption{The $m$-twisted double of a trefoil. In the region marked by
$X$ are $m-3$ full twists.}\lbl{double} 
\end{figure}

\begin{theorem}
\lbl{thm.2}
Fix a knot $K$ in $S^3$ and integers $p,m,n$.
Then $Z_n^{LMO}(\S^p_{D_mK})$ depends only on $p,m$ and the degree
$2n$ part of the Kontsevich integral \cite{Ko} of $K$.
\end{theorem}

With  the above notation, setting $n=1$ we obtain that:
\begin{corollary}
\lbl{cor.n=1}
$$ \l(\S^p_{D_mK})= a_{p,m} \Delta^{''}(K)(1) + b_{p,m}
$$
where $\Delta(K)$ is the Alexander-Conway polynomial of $K$, \cite{C,Ka}
and $a_{p,m}, b_{p,m}$ are constants depending on $p,m$.
\end{corollary}

The above corollary was obtained independently by Ishibe \cite{I}
for general $m,p$, Hoste \cite{Ho} for $m=0$ and Davidow \cite{Da}
for $m=-1, p \equiv \pm 1 \bmod 6$.

We wish to thank Daniel Ruberman for 
numerous encouraging,
enlightening and clarifying conversations.

\section{A reduction of Theorem \ref{thm.1}}
\lbl{sec.skein}

Before we get involved in calculations, we should mention that
the proof of Theorem \ref{thm.1} is an application of  {\em skein theory}
and the following two properties, together with their philosophical proof:
\begin{itemize}
\item[\bf{P1}]
The Casson-Walker-Lescop invariant satisfies a 3-term relation.
This holds since the Casson-Walker-Lescop invariant is a finite type
3 invariant (at least restricted to the set of rational homology 3-spheres,
\cite{GO}).
\item[\bf{P2}]
A crossing change or a smoothening of a link $L$, results in (two) surgeries
along the same  knot in $S^3$. This holds since
the 2-fold branched cover of a disk $D^2$ (branched along
two points) is an annulus, thus the 2-fold branched cover of $D^2 \times I$
branched along two arcs is a solid torus $T$.
\end{itemize}

All links and 3-manifolds in this paper are {\em oriented}.
With an eye in equation \eqref{eq.blah}, we define
for a link $L$ in $S^3$, 
$$ \a(L) \eqdef 1/6 J'_L(-1), \text{\hspace{0.2cm}} \b(L)
\eqdef i^{\s(L)+\n(L)}\l(\S^2_L), \text{\hspace{0.2cm} and \hspace{0.2cm}} 
\g(L) \eqdef  1/4 J_L(-1) \s(L).$$

A triple of links $(L^+, L^-, L^0)$
is called {\em bordered} if there is an embedded disk $D^3$
in $S^3$ that locally intersects them as in figure \ref{crossing}.

\begin{figure}[htpb]
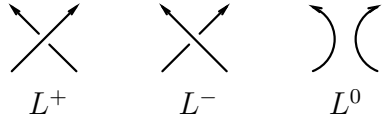

$$ \printname{crossing}
	\setlength{\unitlength}{0.03\standardunitlength}
	\begin{array}{c}  \hspace{-1.7mm}
        	\raisebox{-8pt}{\input draws/crossing.tex }
        	\hspace{-1.9mm}
	\end{array}
 $$
\caption{A bordered triple of links $(L^+, L^-, L^0)$.}\lbl{crossing}
\end{figure}

For a  bordered triple $(L^+, L^-, L^0)$, the skein  
relation $tJ_{L^+}(t) - t^{-1}J_{L^-}(t) =(t^{1/2}-t^{-1/2}) J_{L^0}(t)$
of the Jones polynomial
 implies that:
\begin{equation}
\lbl{eq.a}\begin{split}
\a(L^+) - \a(L^-) = - 2 i \a(L^0) 
+ \frac{J_{L^+}(-1)}{6} + \frac{J_{L^-}(-1)}{6}
\end{split}
\end{equation}

Thus, the following claim:
\begin{claim}
\lbl{claim1}
For a  bordered triple $(L^+, L^-, L^0)$, we have:
\begin{equation}
\begin{split}
\b(L^+)  - \b(L^-) + 2 i \b(L^0)= \g(L^+)  - \g(L^-) + 2 i \g(L^0)   
+ \frac{J_{L^+}(-1)}{6} + \frac{J_{L^-}(-1)}{6}
\end{split}
\end{equation}
\end{claim}
\noindent
together with the initial condition $\a(\text{unknot})=\b(\text{unknot})
-\g(\text{unknot})=0$,
proves Theorem \ref{thm.1}. The rest of the paper is devoted to the proof
of the above claim.

\begin{remark}
\lbl{rem.vj}
Jones and Mullins use a
similar but different skein theory for the Jones polynomial, namely
$t^{-1}V_{L^+}(t) - tV_{L^-}(t) =(t^{1/2}-t^{-1/2}) V_{L^0}(t)$.
The polynomials $V_L$ and $J_L$ are easily seen to be related by:
$J_L(t)=(-1)^{|L|-1}V_L(t^{-1})$, where $|L|$ is the number of components of 
$L$. With our choice, it turns out that
$J_L(1)=2^{|L|-1}$, a positive integer, which is natural from the point of 
view of quantum groups and perturbative Chern-Simons theory.
Furthermore, 
Mullins is evaluating $V_L$ at $-1$ with the the rather nonstandard convention
that $\sqrt{-1}=-i$, whereas we are evaluating $J_L$ at $-1$ with
the convention that $\sqrt{-1}=i$.
\end{remark}

\section{Some linear algebra}

We begin by reviewing three important invariants of symmetric matrices $A$.
All matrices considered have real entries, and $B^T$ denotes the transpose
of the matrix $B$. Two matrices $B$ and $B'$ are called {\em similar}
if $B'=P B P^T$, for a nonsingular matrix $P$. Given a symmetric matrix $A$,
we denote by $\n(A)$, $\s(A)$ its {\em nullity} and {\em signature}
respectively. A lesser known invariant, the {\em sign} $\sgn(A)$
 of $A$,  can be obtained as follows:
bring  $A$ to the form 
$PAP^T=\tByt {A'} 0 0 0 $ 
where $A', P$ are
nonsingular, \cite{Ky}. Then, we can define $\sgn(A)=\sgn(\text{det}(A'))$,
with the understanding that the sign of the determinant of a $0 \times 0$ 
matrix is 1.
It is easy to see that the result is independent of $P$; 
moreover, it coincides with Lescop's
definition \cite[Section 1.3]{Le}. Notice that the signature, nullity
and sign of a matrix do not change under similarity transformations.

We call a  triple of symmetric matrices $(A_+, A_-, A_0)$ {\em bordered} if
$$ A_+ = \tByt {a} {\rho} {\rho^T} {A_0} \text{ and }
A_- = \tByt {a+2} {\rho} {\rho^T} {A_0},
$$ for a row vector $\rho$.
The signatures and nullities  of a bordered triple are related as follows,
\cite{C}:
\begin{equation}
| \n(A_{\pm})-\n(A_0)| + | \s(A_{\pm})-\s(A_0)| =1
\end{equation}
Thus, in a bordered triple, the nullity determines the signature, up to
a sign. A more precise relation is the following: 
\begin{lemma}
\lbl{lem.sign}
The sign, nullity and signature in a bordered triple are related
as follows:
\begin{eqnarray}
\lbl{eq.s}
\s(A_{\pm})-\s(A_0) =
\begin{cases}
0 & \text{ if\hspace{0.2cm} } | \n(A_{\pm})-\n(A_0)|=1 \\ 
\sgn(A_{\pm})\sgn(A_0) & \text{ otherwise. }
\end{cases} \\ \lbl{eq.n}
\n(A_{\pm})-\n(A_0) =
\begin{cases}
0 & \text{ if\hspace{0.2cm} } | \s(A_{\pm})-\s(A_0)|=1 \\
\sgn(A_{\pm})\sgn(A_0) & \text{ otherwise. }
\end{cases}
\end{eqnarray}
Moreover, if $\e_x \eqdef \sgn(A_x) i^{\s(A_x)+\n(A_x)}$ for
$x \in \{ +, -, 0 \}$, then we have:
\begin{equation}
\lbl{eq.eta}
\e_+=\e_-= i \e_0.
\end{equation}
\end{lemma}

\begin{proof}
By similarity transformations, we can assume that:
$$
A_+= \tByt {a} {\rho} {\rho^T} {0} \oplus D,
A_+= \tByt {a+2} {\rho} {\rho^T} {0} \oplus D, A_0= D \oplus [0]^r
$$ 
where $D$ is a nonsingular, diagonal matrix, $[0]^r$ is the zero $r \times
r$ matrix,
$\rho$ is a $1\times r$ vector and $a$ a real number.
Since  the nullity, signature and sign
of the matrix  $\tByt {a} {\rho} {\rho^T} {0}$ are given by: 
\begin{center}
\begin{tabular}{|l|c|c|c|}
\hline
        & $\rho=a=0$ & $\rho=0, a \neq 0$ & $ \rho \neq 0$ \\ \hline
nullity & $r+1$      &     $r$            &  $r-1$         \\ \hline
signature &  $0$     &   $\sgn(a)$        &   $0$           \\ \hline
sign &       $1$     &   $\sgn(a)$        &   $-1$          \\ \hline
\end{tabular}
\end{center}
\noindent
the result follows by a case-by-case argument.
\end{proof}

\begin{remark}
\lbl{rem.det}
For future reference, we mention that the determinants of a bordered triple
of matrices are related by:
\begin{equation}
\lbl{eq.det}
 \det(A_+) - \det(A_-) + 2 \det(A_0) =0.
\end{equation}
This follows easily by expanding the first two determinants along the
first column.
\end{remark}

Given an oriented link $L$ in $S^3$, choose of Seifert surface of it, together
with a basis for its homology and and consider the associated Seifert matrix
$E_L$. Recall that the {\em nullity} $\n(L)$, {\em signature} 
$\s(L)$ and sign $\sgn(E_L)$ of $L$ are defined as the nullity, signature
and sign of the symmetrized Seifert matrix $E_L + E_L^T$. It turns out that
the signature and nullity of a link are independent of the Seifert surface
chosen, and that $\n(L)=\b_1(\S^2_L)$, where $\b_1$ is the first betti number.
On the other hand,  $\sgn(E_L)$ depends on the Seifert matrix.  

It is easy to see that given a bordered triple $(L^+, L^-, L^0)$ of links, 
one can construct a triple of 
Seifert matrices so that the associated triple of symmetrized Seifert
matrices is bordered, \cite{C}.

\section{Proof of Theorem \ref{thm.1}}
\lbl{sec.lescop}
 
\subsection{The Casson-Walker-Lescop invariant of 3-manifolds}
\lbl{sub.gen}

Given an (integrally) framed oriented r-component link $\L$ in $S^3$ 
(with ordered components), 
let $S^3_{\L}$ denote the closed 3-manifold obtained by Dehn surgery
on $\L$. Its linking matrix, $F(\L)$ gives a presentation
of $H_1(S^3_{\L}, \BZ)$. Notice that $\n(F(\L))= \b_1(S^3_{\L})$.
The Casson-Walker-Lescop invariant $\l$ of $S^3_{\L}$ is defined by:
$$
\l(S^3_{\L})=\sgn(F(\L))(D(\L)+ H_0(\L)+ H_1(\L) + H_2(\L)), \text{ where } 
$$
\begin{eqnarray*}
D(\L) & = & \sum_{\emptyset \neq \L' \subseteq \L} 
            \det(F(\L \setminus \L')) \zeta(\L'), \\
H_0(\L) & = & \frac{\det(F(\L))}{4}\s(F(\L)) \\
H_1(\L) & = & - \frac{1}{6} \sum_{j=1}^r \det(F(\L \setminus {j})) \\
H_2(\L) & = &  \frac{1}{12} \sum_{\emptyset \neq \L' \subseteq \L}   
             \det(F(\L \setminus \L'))  (-1)^{|\L'|} L_8(\L'),
\end{eqnarray*}
$\zeta(\L)$ is a special value of (a derivative of) the multivariable
Alexander polynomial of $\L$ and  $L_8(\L)$ is a polynomial in the linking
numbers $l_{ab}$ ($a,b=1 \dots r$) of $\L$ given explicitly by:
$$
L_8(\L)=\sum_{j=1}^r\sum_{\s \in \text{Sym}_r}
l_{i\s(1)} l_{\s(1)\s(2)} \dots l_{\s(r-1)\s(r)}l_{\s(r)i},
$$
where $\text{Sym}_r$ is the symmentric group with $r$ letters.

Notice also that since the links $\L$ that we consider will be
 integrally framed, the Dedekind sums
appearing in \cite[definition 1.4.5]{Le} vanish.

\subsection{A construction of 2-fold branched covers}
\lbl{sub.branched}

In this section we review the details of some well known construction
of 2-fold branched covers of links in $S^3$. For a general reference,
see \cite{Ka,AK}. Given an (oriented) link $L$ in $S^3$,
choose a Seifert surface $F_L$ of $L$, and a basis of its first
homology  and let $E_L$ be its Seifert matrix.
Push a bicolar of  $F_L$ in the interior of $D^4$
(the 4-manifold obtained is still diffeomorphic to $D^4$), and glue
two copies of the obtained 4-manifold along $F_L$ according to the pattern of
\cite[p. 281]{Ka}; let $N_{F_L}$ denote the resulting 4-manifold, which
is a 2-fold cover of $D^4$ branched along $F_L$. Its boundary is
$\S^2_L$, the 2-fold cover of $S^3$ branched along $L$.
In \cite[section 2]{AK}, Akbulut-Kirby showed that $N_{F_L}$ is a 
4-dimensional handlebody (i.e. the result of attaching 2-handles along $D^4$),
and that the intersection form with respect to some basis of these
2-handles is the symmetrized Seifert matrix $E_L + E_L^T$ of $L$.
Let $\L$ denote the cores in $S^3$ of the 2-handles. Thus, $\L$ is
a framed link in $S^3$ with linking matrix $E_L + E_L^T$, such that
Dehn surgery on $\L$ is $\S^2_L$. Of course, the link $\L$ depends on 
the choice of Seifert surface of $L$ as well as on a choice of basis
on its homology. Akbulut-Kirby \cite{AK} describe an algorithm for
drawing $\L$ and implement it with beautiful pictures, however we will not
need the precise picture of the link $\L$!

Assume now that $(L^+, L^-, L^0)$ is a bordered triple with
admissible Seifert surfaces. 
Property {\bf P1} of Section \ref{sec.skein} implies that
there is a solid torus in
$\S^2_{L^0}$ and three simple curves $\a^+, \a^-, \a^0$ in its boundary
so that $\S^2_{L^x}$ is diffeomorphic to surgery on the solid torus $T$
in $\S^2_{L^0}$ along $\a_x$, for $x \in \{+,-,0 \}$. Using the argument
of \cite[page p.429]{Mu}, it follows that there is
a choice of a standard symplectic basis $\{ x_1, x_2 \}$ for $H_1(\partial T,
\BZ)$ so that
$\a^+= x_2, \a^-= 2 x_1 + x_2,  \a^0= x_1$. In other words, we have:
$\la \a^0, \a^+ \ra =1$ and $\a^-= 2 \a^0 + \a^+$, where $\la \cdot, \cdot \ra$
is the intersection form. 

From the above discussion, it follows that there is a triple
 of framed 
oriented links (not necessarily bordered!) 
$(\L^+, \L^-, \L^0)$ in $S^3$ an oriented knot $K$ and an  integer $n$
 so that $\S^2_{L^x}=\S^3_{\L^x}$ for all $x \in \{ +, -, 0 \}$,
and so that $\L^+$ (resp. $\L^-$)
 is the disjoint union of $\L^0$ with the framed knot
$(K,n)$ (resp. $(K, n+2)$). Thus, the triple  
$(F(\L^+), F(\L^-), F(\L^0))$ of linking matrices
is bordered, with nullity, signature and sign equal to that of
the triple $(L^+, L^-, L^0)$.

\subsection{Proof of Claim \ref{claim1}}

Let us define $\e_x=\sgn(\L^x)i^{\s(L^x)+\n(L^x)}$ for $x \in \{ -, +, 0 \}$, 
Since the $\s(\L^x)=\s(L^x)$,
Lemma \ref{lem.sign} implies that $\e_{\pm}= i \e_0$.

We can now calculate as follows:
\begin{eqnarray*}
\b(L^+)-\b(L^-)+2i\b(L^0) & = &
 \big\{ \e_+ D(\L^+) - \e_- D(\L^-) +2i \e_0 D(\L^0) \big\} \\ 
&  + &\sum_{k=0}^2
 \big\{   \e_+ H_k(\L^+) - \e_- H_k(\L^-) +2i \e_0 H_k(\L^0) \big\}
\end{eqnarray*}
\noindent
We  claim that:
\begin{eqnarray}
\lbl{eq.1}
 \e_+ D(\L^+) - \e_- D(\L^-) +2i \e_0 D(\L^0) & = & 0 \\ \lbl{eq.2}
\e_+ H_0(\L^+) - \e_- H_0(\L^-) +2i \e_0 H_0(\L^0) & = & 
\g(L^+)  - \g(L^-) +2 i \g(L^+) \\ \lbl{eq.3}
\e_+ H_1(\L^+) - \e_- H_1(\L^-) +2 \e_0 H_1(\L^0) & = & 0 \\ \lbl{eq.4}
\e_+ H_2(\L^+) - \e_- H_2(\L^-) +2 \e_0 H_2(\L^0) & = & 
 \frac{J_{L^+}(-1)}{6} + \frac{J_{L^-}(-1)}{6}
\end{eqnarray}

Before we show the above equations, we let 
$\L^0 \cup K^{\pm}= \L^{\pm}$, and $I$ denote an arbitrary nonempty
sublink of $\L^0$ with complement $I' \eqdef \L^0 \setminus I$. 

Using equation \eqref{eq.eta}, it follows that
the left hand side of \eqref{eq.1} equals to:
\begin{eqnarray*}
 \sum_{I}  & & \big\{ \e_+ \det(F(I' \cup K^+)) - \e_-
 \det(F(I' \cup K^-)) +2 i \e_0 \det(F(I')) \big\} \zeta(I) + \\ 
 \sum_{I}  & &
\big\{ \e_+  \det(F(I')) - \e_-
 \det(F(I')) \big) \big\} \zeta(I \cup K) 
\end{eqnarray*}
Using equation \eqref{eq.eta}, 
and the fact that $F(I \cup K^+), F(I \cup K^-), F(I \cup K^0))$
is a bordered triple of matrices, it follows by Remark \ref{rem.det}
that the first and second sum shown above vanishes, 
thus showing equation \eqref{eq.1}.

In order to show equation \eqref{eq.2}, use the fact
 that for a link $L$ in $S^3$
we have: 
$$|J_L(-1)|=|H_1(\S^2, \BZ)|= i^{-\s(L)- 2 \n(L)}J_L(-1)$$
 (compare with
\cite[Theorem 2.4]{Mu},
and with Kauffman \cite{Ka}, with the understanding that the order
of an infinite group is $0$, and keeping in mind Remark \ref{rem.vj}). Thus,
\begin{equation}
\lbl{eq.h1}
\e_x \det(F(\L^x)) = i^{\s(L^x) + \n(L^x)} |H_1(\S^2, \BZ)|=
i^{-\n(L^x)} J_{L^x}(-1) = J_{L^x}(-1)
\end{equation}
 for all $x \in \{ +, -, 0 \}$ (where the last equality above follows
from the fact that if $\n(L^x) \neq 0$, then $J_{L^x}(-1)=0$).   
Since $\s(\L^x)=\s(L^x)$ for all $x$, it follows 
that $\e_x H_0(\L^x) = \g(L^x)$ for all $x$, which proves equation 
\eqref{eq.2}.

Equation \eqref{eq.3} follows in the same way as equation \eqref{eq.1}
shown above.

Using \eqref{eq.eta} and the definition of $H_2$, it follows that
the left hand side of equation \eqref{eq.4} equals to:
\begin{eqnarray*}
& i \e_0  & (H_2(\L^+) - H_2(\L^0) + 2 H_2(\L^0)) \\ 
= & \frac{i \e_0}{12} & \big\{ 
\sum_I \det(F(I')) (-1)^{|I\cup K|} L_8(I \cup K^+)
+ \sum_I \det(F(I' \cup K^+)) (-1)^{|I|} L_8(I) \\
& - &
\sum_I \det(F(I')) (-1)^{|I\cup K|} L_8(I \cup K^-)
- \sum_I \det(F(I' \cup K^-)) (-1)^{|I|} L_8(I) \\
& + &
2 \sum_I \det(F(I')) (-1)^{|I|} L_8(I) \big\} \\
= & \frac{i \e_0}{12} & \big\{
\sum_I \det(F(I')) (-1)^{|I|} ( -  L_8(I \cup K^+) + L_8(I \cup K^+))
\big\}
\end{eqnarray*}
\noindent
It is easy to see that if $(l^{\pm}_{ab})=F(\L^{\pm})$ and that $K^{\pm}$
is the first ordered component of $\L^{\pm}$, then
\begin{eqnarray*}
-  L_8(I \cup K^+) + L_8(I \cup K^+))
& 
= &
2 \sum_{\s''}
l^+_{1\s''(1)}l^+_{\s''(1)\s''(2)} \dots l^+_{\s''(r-1)\s''(r)} \\ 
& + & 2 \sum_{\s'} 
l^-_{\s'(1)\s'(2)} l^-_{\s'(2)\s'(3)} \dots l^-_{\s'(r+1)1}  
\end{eqnarray*}
where the summation is over all $\s',  \in \text{Sym}_{r+1}$ (resp. $\s''$)
such that $\s'(1)=1$ (resp. $\s''(r+1)=1$). Combined with the above, 
and with \eqref{eq.eta}, \eqref{eq.h1}, 
 the left hand side of equation \eqref{eq.4} equals to
\begin{eqnarray*}
\frac{i \e_0}{6}(\det(F(\L^+))+\det(F(\L^-))=
\frac{1}{6}(\e_+ \det(F(\L^+)) + \e_- \det(F(\L^-))=
\frac{1}{6}(V_{L^+}(-1)+ V_{L-}(-1))
\end{eqnarray*} 
which concludes the proof of equation \eqref{eq.4} and of Claim
\ref{claim1}.


Corollary \ref{cor.1} follows from the fact that the Casson-Walker-Lescop
invariant of a manifold with betti number at least 4, vanishes, \cite{Le}.
Remark \ref{rem.mullins} follows from Theorem \ref{thm.1} and
\eqref{eq.h1}.

\begin{remark}
The forth roots of unity in \eqref{eq.eta} work out in such a way
as to obtain 
the cancellations in equations \eqref{eq.1}, \eqref{eq.2}, \eqref{eq.3}
and \eqref{eq.4}.  
The philosophical reason for the cancelation in
equation \eqref{eq.1} is property {\bf P2} of Section \ref{sec.skein}.
\end{remark}

\section{Proof of Theorem \ref{thm.2}}

Fix a knot $K$ in $S^3$ and integers $p,m,n$. Our first goal is to give a 
Dehn surgery description of the $p$-fold branched cover 
$\S^p_{D_mK}$ of the $m$-twisted double $D_mK$ of $K$.

\begin{figure}[htpb]
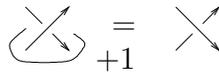

$$ \printname{Dehn}
	\setlength{\unitlength}{0.03\standardunitlength}
	\begin{array}{c}  \hspace{-1.7mm}
        	\raisebox{-8pt}{\input draws/Dehn.tex }
        	\hspace{-1.9mm}
	\end{array}
 $$
\caption{A Dehn twist along a $+1$ framed unknot 
intersecting two arcs (in an arbitrary
3-manifold) on the left
hand side gives a diffeomorphic 3-manifold with a different
embedding of the two arcs.}\lbl{Dehn}
\end{figure}

We begin with a definition: a link $L$ in $S^3$ is called $K$-{\em unknotting}
if it is  unit-framed,  algebraically split 
(i.e., one with linking numbers zero and framing $\pm 1$ on each component),
lying in a standard solid torus in $S^3$, 
such that Dehn surgery $S^3_L$ on $L$ is diffeomorphic to $S^3$ and such
that the image of a meridian of the solid torus in $S^3_L$ is isotopic to
the knot $K$ in $S^3$. 
A $K$-uknotting link $L$ can be obtained by projecting $K$ to 
a generic plane, choosing a set of crossings that unknot $K$,
and placing an unknot with framing $\pm 1$ around each crossing. The union
of these unknots is the desired link $L$, as follows from Figure \ref{Dehn}.

We represent a $K$-unknotting link $L$ 
in the standard solid torus in $S^3$ with
the left hand side of Figure \ref{Dehn}. On the right hand side of the
same figure is shown the $m$-twisted double of $K$.

\begin{figure}[htpb]
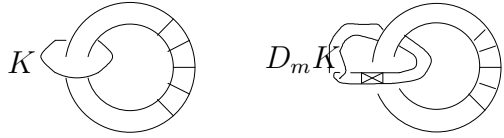

$$ \printname{step1}
	\setlength{\unitlength}{0.03\standardunitlength}
	\begin{array}{c}  \hspace{-1.7mm}
        	\raisebox{-8pt}{\input draws/step1.tex }
        	\hspace{-1.9mm}
	\end{array}
 $$
\caption{On the left, a $K$-unknotting link $L$.
On the right, the  $m$-twisted double $D_mK$. In the box marked $X$ are
$m$ twists.}
\end{figure}

Next, we construct an $D_mK$ uknotting link:
using Figure \ref{Dehn}, we introduce an extra unit-framed unknot $C$
to unknot $D_mK$ as shown on the left hand side of Figure \ref{step2}, 
and we isotope the result as shown on the right hand side of the figure.
Then, $L \cup C$ is a $D_mK$-unknotting link.

\begin{figure}[htpb]
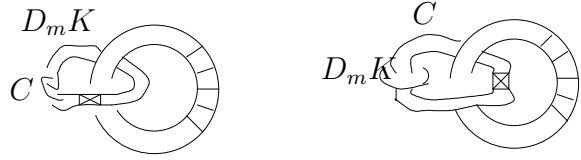

$$ \printname{step2}
	\setlength{\unitlength}{0.03\standardunitlength}
	\begin{array}{c}  \hspace{-1.7mm}
        	\raisebox{-8pt}{\input draws/step2.tex }
        	\hspace{-1.9mm}
	\end{array}
 $$
\caption{Two isotopic views of $L \cup C$,  which is  
$D_mK$-unknotting.}\lbl{step2}
\end{figure}

Cutting the meridian disk that $D_mK$ bounds on the right hand side of
Figure \ref{step2}, and gluing $p$ of them side by side as in Figure
\ref{pfold}, gives a framed link $L(p,m)$ in $S^3$, which is a
Dehn surgery presentation of $\S^p_{D_mK}$. 

\begin{figure}[htpb]
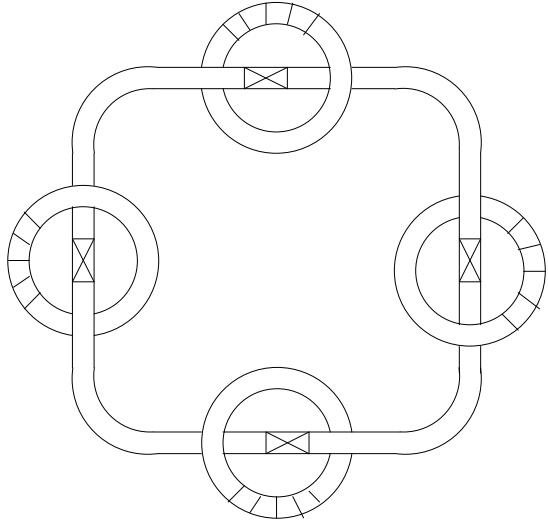

$$ \printname{step4}
	\setlength{\unitlength}{0.03\standardunitlength}
	\begin{array}{c}  \hspace{-1.7mm}
        	\raisebox{-8pt}{\input draws/step4.tex }
        	\hspace{-1.9mm}
	\end{array}
 $$  
\caption{The framed link $L(4,m)$.}\lbl{pfold}
\end{figure}

Next, we review some well-known facts about the LMO invariant, and its
cousin, the Aarhus integral, \cite{BGRT}. We will assume some familiarity
with the standard definitions of finite type invariants of links and
3-manifolds, \cite{B-N,Ko,LMO,BGRT}. Both the LMO invariant and
the Aarhus integral of a 3-manifold obtained by Dehn surgery of a link
in $S^3$ are defined in terms of the Kontsevich integral of the link.
For rational homology 3-spheres, the LMO invariant equals, properly
normalized, to the Aarhus integral, where the normalization of the factor
is the order of the first homology group with integer coefficients. 
The manifolds in question, namely $\S^p_{D_mK}$, are rational homology
3-spheres, and the order of the first homology group depends only on
$m$ and $p$ (since the first homology is given in terms of the Alexander
polynomial of $D_mK$, which depends only on $m$, evaluated at $p$th
roots of unity).
Thus, it suffices to show that the degree at most $n$ 
part of the Aarhus integral
$Z_n^A(\S^p_{D_mK})$ of $\S^p_{D_mK}$ depends on the degree at most $2n$ 
part of the Kontsevich integral $Z_{2n}^K(K)$ of $K$.

Since the degree at most $d$ part of the Kontsevich integral is the universal
Vassiliev invariant of type $d$ \cite{Ko, B-N},  it suffices to show that
the knot invariant $ K \to Z^{A}_n(\S^p_{D_mK})$ is a Vassiliev
invariant of type $2n$.

A crossing change on $K$ corresponds to an extra
component to the $K$-unknotting link $L$ and $p$ extra components to the 
link $L(p,m)$, thus the alternating sum of 
$2n+1$ double points on $K$, corresponds to 
 the alternating sum of a link $L_{\text{alt}} \cup L(p,m)$, where
$L_{\text{alt}}$ has $(2n+1)p$ components, 
and where we alternate by considering
including or not each of the $p$ components. We then consider the Kontsevich
integral of $L_{\text{alt}} \cup L(p,m)$, expressed in terms of disjoint union
of uni-trivalent graphs, the legs of which we glue according to the
definition of the Aarhus integral. Since the linking numbers between 
the components of $L_{\text{alt}}$ and $L_{\text{alt}} \cup L(p,m)$ are zero, 
a standard counting
argument (compare with \cite{L}, \cite[part II, Section 4, proof of 
Theorem 2]{BGRT} 
and with
\cite[proof of Theorem 3]{GH})shows that after alernating $2n+1$ terms, 
the degree $n$ part of $Z^A_n(L_{\text{alt}} \cup L(p,m))$ vanishes.

Corollary \ref{cor.n=1} follows immediately from Theorem \ref{thm.2},
using the fact that the degree $1$ part of the LMO invariant equals to
the Casson-Walker-Lescop invariant, \cite{LMMO,LMO}.

\ifx\undefined\bysame
	\newcommand{\bysame}{\leavevmode\hbox to3em{\hrulefill}\,}
\fi

\end{document}

%% file: draws/double.tex
\begingroup\makeatletter\ifx\SetFigFont\undefined%
\gdef\SetFigFont#1#2#3#4#5{%
  \reset@font\fontsize{#1}{#2pt}%
  \fontfamily{#3}\fontseries{#4}\fontshape{#5}%
  \selectfont}%
\fi\endgroup%
{\renewcommand{\dashlinestretch}{30}
\begin{picture}(3380,3360)(0,-10)
\put(1315.500,3170.500){\arc{325.960}{4.7892}{7.2873}}
\put(1340.500,3020.500){\arc{325.960}{2.1375}{4.6356}}
\path(3053,2433)(3353,2433)
\path(3128,2133)(3353,2133)
\path(3053,2433)(3353,2133)
\path(3353,2433)(3128,2133)
\path(653,1533)(655,1534)(659,1537)
	(667,1542)(678,1550)(694,1560)
	(715,1574)(740,1590)(769,1608)
	(802,1629)(836,1651)(873,1673)
	(910,1695)(948,1717)(985,1738)
	(1021,1758)(1055,1775)(1088,1791)
	(1119,1804)(1148,1815)(1176,1823)
	(1202,1829)(1228,1832)(1253,1833)
	(1278,1832)(1303,1829)(1328,1823)
	(1353,1816)(1378,1806)(1403,1794)
	(1428,1781)(1453,1766)(1478,1750)
	(1503,1733)(1528,1715)(1553,1697)
	(1578,1679)(1603,1660)(1628,1643)
	(1653,1626)(1678,1609)(1703,1594)
	(1728,1581)(1753,1568)(1778,1557)
	(1803,1548)(1828,1540)(1853,1533)
	(1878,1527)(1904,1523)(1930,1519)
	(1958,1516)(1987,1514)(2018,1513)
	(2051,1512)(2085,1513)(2121,1513)
	(2158,1515)(2196,1516)(2233,1518)
	(2270,1520)(2304,1522)(2337,1524)
	(2366,1526)(2391,1528)(2412,1530)
	(2428,1531)(2439,1532)(2447,1533)
	(2451,1533)(2453,1533)
\path(653,1233)(655,1234)(659,1237)
	(667,1242)(678,1250)(694,1260)
	(715,1274)(740,1290)(769,1308)
	(802,1329)(836,1351)(873,1373)
	(910,1395)(948,1417)(985,1438)
	(1021,1458)(1055,1475)(1088,1491)
	(1119,1504)(1148,1515)(1176,1523)
	(1202,1529)(1228,1532)(1253,1533)
	(1278,1532)(1303,1529)(1328,1523)
	(1353,1516)(1378,1506)(1403,1494)
	(1428,1481)(1453,1466)(1478,1450)
	(1503,1433)(1528,1415)(1553,1397)
	(1578,1379)(1603,1360)(1628,1343)
	(1653,1326)(1678,1309)(1703,1294)
	(1728,1281)(1753,1268)(1778,1257)
	(1803,1248)(1828,1240)(1853,1233)
	(1878,1227)(1904,1223)(1930,1219)
	(1958,1216)(1987,1214)(2018,1213)
	(2051,1212)(2085,1213)(2121,1213)
	(2158,1215)(2196,1216)(2233,1218)
	(2270,1220)(2304,1222)(2337,1224)
	(2366,1226)(2391,1228)(2412,1230)
	(2428,1231)(2439,1232)(2447,1233)
	(2451,1233)(2453,1233)
\path(1628,1308)(1627,1306)(1624,1302)
	(1618,1295)(1610,1285)(1598,1270)
	(1584,1252)(1567,1230)(1547,1206)
	(1525,1179)(1502,1152)(1478,1124)
	(1454,1096)(1429,1069)(1406,1045)
	(1383,1022)(1360,1001)(1338,983)
	(1317,967)(1296,954)(1275,942)
	(1253,933)(1232,926)(1211,920)
	(1189,915)(1166,911)(1142,908)
	(1118,905)(1093,903)(1068,902)
	(1042,901)(1016,901)(990,901)
	(964,901)(937,902)(911,903)
	(884,904)(858,906)(831,907)
	(805,909)(779,912)(754,915)
	(728,918)(703,922)(678,927)
	(653,933)(628,940)(603,947)
	(577,955)(551,964)(524,973)
	(497,982)(469,991)(441,1001)
	(412,1011)(383,1020)(354,1030)
	(325,1041)(296,1052)(268,1063)
	(240,1075)(213,1087)(188,1101)
	(163,1115)(141,1131)(119,1148)
	(100,1167)(82,1187)(67,1209)
	(53,1233)(44,1254)(36,1276)
	(29,1301)(23,1327)(19,1356)
	(16,1388)(14,1422)(12,1459)
	(12,1499)(12,1541)(13,1585)
	(15,1630)(17,1677)(20,1725)
	(23,1773)(26,1820)(30,1866)
	(33,1910)(37,1951)(40,1989)
	(43,2022)(45,2052)(48,2076)
	(49,2095)(51,2110)(52,2121)
	(52,2128)(53,2131)(53,2133)
\path(1403,1533)(1402,1530)(1400,1525)
	(1397,1516)(1391,1501)(1384,1483)
	(1376,1460)(1366,1435)(1356,1408)
	(1344,1381)(1333,1354)(1321,1329)
	(1310,1306)(1299,1286)(1288,1268)
	(1277,1254)(1265,1242)(1253,1233)
	(1238,1225)(1222,1219)(1202,1215)
	(1180,1213)(1156,1212)(1129,1213)
	(1100,1215)(1071,1218)(1042,1221)
	(1015,1224)(993,1227)(975,1230)
	(963,1231)(956,1233)(953,1233)
\path(953,1233)(951,1233)(945,1232)
	(936,1232)(922,1230)(902,1229)
	(878,1227)(848,1225)(815,1223)
	(779,1220)(742,1218)(703,1216)
	(666,1214)(629,1213)(594,1212)
	(562,1213)(533,1213)(507,1215)
	(483,1218)(462,1222)(444,1227)
	(428,1233)(410,1242)(395,1253)
	(381,1265)(369,1279)(358,1293)
	(349,1309)(340,1325)(332,1341)
	(325,1358)(318,1376)(312,1394)
	(305,1412)(300,1431)(294,1450)
	(289,1469)(284,1490)(281,1511)
	(278,1533)(277,1552)(277,1573)
	(277,1596)(279,1621)(282,1648)
	(286,1678)(290,1710)(295,1744)
	(301,1780)(308,1817)(314,1854)
	(321,1891)(327,1926)(333,1957)
	(339,1986)(343,2009)(347,2028)
	(350,2042)(352,2051)(353,2056)(353,2058)
\path(1703,1608)(1704,1611)(1707,1616)
	(1712,1626)(1719,1640)(1728,1658)
	(1739,1679)(1751,1702)(1763,1725)
	(1776,1747)(1789,1768)(1802,1786)
	(1814,1802)(1827,1815)(1839,1825)
	(1853,1833)(1868,1839)(1883,1844)
	(1900,1848)(1918,1851)(1937,1853)
	(1956,1855)(1976,1857)(1996,1858)
	(2016,1858)(2036,1858)(2057,1857)
	(2076,1855)(2096,1852)(2115,1847)
	(2134,1841)(2153,1833)(2168,1825)
	(2183,1815)(2199,1805)(2215,1793)
	(2231,1781)(2248,1768)(2265,1754)
	(2282,1740)(2300,1726)(2317,1711)
	(2334,1696)(2351,1680)(2367,1664)
	(2383,1648)(2398,1631)(2411,1613)
	(2424,1595)(2435,1575)(2445,1555)
	(2453,1533)(2459,1514)(2463,1493)
	(2467,1470)(2470,1445)(2472,1418)
	(2473,1388)(2474,1356)(2473,1321)
	(2473,1285)(2471,1246)(2470,1207)
	(2468,1168)(2466,1129)(2464,1092)
	(2462,1058)(2460,1027)(2458,1000)
	(2456,978)(2455,960)(2454,948)
	(2453,940)(2453,935)(2453,933)
\path(1478,1833)(1479,1836)(1482,1841)
	(1487,1850)(1494,1864)(1503,1881)
	(1514,1901)(1526,1923)(1538,1945)
	(1551,1967)(1564,1988)(1577,2006)
	(1589,2022)(1602,2036)(1614,2048)
	(1628,2058)(1643,2067)(1659,2075)
	(1676,2081)(1694,2087)(1713,2093)
	(1733,2098)(1754,2102)(1774,2106)
	(1795,2110)(1816,2114)(1837,2117)
	(1856,2121)(1876,2124)(1894,2127)
	(1911,2130)(1928,2133)(1946,2135)
	(1964,2137)(1980,2139)(1995,2140)
	(2010,2141)(2023,2143)(2037,2144)
	(2051,2144)(2065,2144)(2080,2144)
	(2096,2143)(2113,2141)(2132,2138)
	(2153,2133)(2168,2129)(2185,2124)
	(2202,2118)(2220,2112)(2239,2105)
	(2259,2098)(2280,2090)(2301,2082)
	(2322,2074)(2344,2065)(2366,2056)
	(2387,2048)(2408,2039)(2428,2030)
	(2448,2022)(2466,2014)(2483,2006)
	(2499,1998)(2514,1991)(2528,1983)
	(2546,1972)(2561,1962)(2575,1951)
	(2587,1941)(2598,1931)(2608,1921)
	(2617,1912)(2626,1902)(2635,1891)
	(2643,1881)(2652,1870)(2661,1858)
	(2670,1846)(2678,1833)(2686,1819)
	(2693,1806)(2699,1793)(2706,1780)
	(2711,1768)(2717,1757)(2723,1745)
	(2728,1733)(2733,1719)(2738,1703)
	(2742,1685)(2746,1663)(2750,1638)
	(2753,1608)(2754,1588)(2756,1567)
	(2757,1542)(2757,1515)(2758,1485)
	(2758,1452)(2758,1416)(2758,1377)
	(2758,1336)(2758,1292)(2757,1248)
	(2757,1203)(2756,1159)(2756,1116)
	(2755,1077)(2755,1041)(2754,1010)
	(2754,984)(2754,964)(2753,950)
	(2753,941)(2753,935)(2753,933)
\path(653,933)(653,930)(654,924)
	(654,913)(656,897)(657,875)
	(660,849)(662,820)(666,788)
	(669,756)(674,724)(679,693)
	(684,665)(691,639)(698,615)
	(706,594)(716,575)(728,558)
	(740,543)(753,529)(767,516)
	(781,502)(795,490)(810,477)
	(824,464)(839,452)(854,440)
	(870,428)(886,416)(903,405)
	(921,394)(941,383)(962,373)
	(985,363)(1010,354)(1038,346)
	(1069,339)(1103,333)(1126,330)
	(1150,327)(1176,324)(1204,322)
	(1233,320)(1263,318)(1295,317)
	(1329,315)(1363,314)(1399,313)
	(1436,312)(1474,311)(1513,310)
	(1552,309)(1592,309)(1632,308)
	(1671,308)(1711,308)(1750,308)
	(1788,308)(1826,308)(1863,309)
	(1899,310)(1933,311)(1966,312)
	(1998,314)(2028,316)(2056,319)
	(2083,321)(2108,325)(2131,329)
	(2153,333)(2185,341)(2213,351)
	(2239,362)(2261,375)(2281,388)
	(2298,403)(2314,419)(2328,435)
	(2340,452)(2352,469)(2363,486)
	(2374,504)(2385,521)(2395,538)
	(2405,555)(2416,571)(2426,587)
	(2436,603)(2445,618)(2453,633)
	(2461,652)(2467,672)(2471,693)
	(2473,716)(2474,740)(2473,767)
	(2471,795)(2468,823)(2465,850)
	(2462,875)(2459,896)(2456,912)
	(2455,924)(2453,930)(2453,933)
\path(353,1008)(353,1005)(352,999)
	(352,989)(351,973)(350,952)
	(348,926)(347,898)(346,867)
	(344,835)(343,803)(343,773)
	(343,745)(343,718)(345,694)
	(346,672)(349,652)(353,633)
	(358,614)(363,596)(368,579)
	(374,562)(381,545)(387,528)
	(393,512)(400,496)(407,479)
	(414,463)(422,447)(430,431)
	(439,415)(449,399)(461,383)
	(473,366)(487,350)(503,333)
	(517,319)(533,306)(550,291)
	(568,277)(587,262)(608,247)
	(629,232)(651,216)(674,201)
	(697,185)(721,169)(745,153)
	(768,138)(791,123)(814,109)
	(836,95)(858,82)(879,70)
	(899,59)(917,49)(936,40)
	(953,33)(976,25)(999,19)
	(1022,15)(1046,13)(1071,12)
	(1097,13)(1124,15)(1151,18)
	(1176,21)(1199,24)(1219,27)
	(1234,30)(1244,31)(1250,33)(1253,33)
\path(1253,33)(1255,33)(1258,33)
	(1265,33)(1275,33)(1290,33)
	(1309,32)(1333,32)(1362,32)
	(1395,31)(1433,31)(1473,31)
	(1517,30)(1562,30)(1609,29)
	(1657,29)(1704,29)(1751,28)
	(1796,28)(1840,28)(1882,28)
	(1922,28)(1959,28)(1993,28)
	(2025,29)(2055,29)(2082,30)
	(2108,31)(2131,32)(2153,33)
	(2187,35)(2218,38)(2246,41)
	(2271,44)(2294,47)(2316,50)
	(2336,53)(2354,57)(2372,60)
	(2390,64)(2407,67)(2424,72)
	(2441,76)(2458,81)(2476,87)
	(2493,93)(2511,100)(2528,108)
	(2546,118)(2564,128)(2580,139)
	(2596,150)(2612,160)(2628,171)
	(2643,181)(2658,192)(2673,204)
	(2687,216)(2700,230)(2713,246)
	(2725,263)(2736,283)(2745,307)
	(2753,333)(2757,353)(2761,374)
	(2764,399)(2766,426)(2767,456)
	(2768,489)(2768,525)(2768,564)
	(2768,605)(2767,649)(2766,693)
	(2764,738)(2763,782)(2761,825)
	(2759,864)(2758,900)(2757,931)
	(2755,957)(2755,977)(2754,991)
	(2753,1000)(2753,1006)(2753,1008)
\path(2753,1608)(2756,1609)(2761,1612)
	(2770,1616)(2785,1623)(2803,1632)
	(2826,1644)(2851,1657)(2878,1672)
	(2905,1687)(2932,1704)(2957,1720)
	(2980,1737)(3000,1755)(3018,1772)
	(3032,1791)(3044,1811)(3053,1833)
	(3059,1851)(3063,1871)(3067,1892)
	(3070,1917)(3072,1944)(3073,1973)
	(3074,2005)(3073,2040)(3073,2077)
	(3071,2115)(3070,2154)(3068,2194)
	(3066,2233)(3064,2271)(3062,2306)
	(3060,2338)(3058,2365)(3056,2387)
	(3055,2405)(3054,2418)(3053,2426)
	(3053,2431)(3053,2433)
\path(2753,1308)(2756,1308)(2762,1309)
	(2773,1310)(2789,1312)(2811,1314)
	(2837,1317)(2866,1320)(2898,1324)
	(2930,1329)(2962,1333)(2993,1339)
	(3021,1344)(3047,1350)(3071,1357)
	(3092,1365)(3111,1373)(3128,1383)
	(3144,1394)(3160,1406)(3175,1419)
	(3189,1432)(3204,1445)(3218,1459)
	(3231,1473)(3245,1487)(3258,1501)
	(3271,1516)(3284,1532)(3296,1549)
	(3308,1567)(3319,1586)(3329,1607)
	(3338,1630)(3346,1655)(3353,1683)
	(3357,1704)(3360,1727)(3363,1752)
	(3365,1779)(3367,1809)(3368,1842)
	(3368,1877)(3368,1916)(3368,1957)
	(3368,2000)(3367,2044)(3366,2090)
	(3364,2135)(3363,2181)(3361,2224)
	(3360,2265)(3358,2302)(3357,2336)
	(3356,2364)(3355,2387)(3354,2405)
	(3354,2418)(3353,2426)(3353,2431)(3353,2433)
\path(1328,3333)(1326,3333)(1322,3333)
	(1314,3333)(1302,3333)(1286,3333)
	(1265,3332)(1239,3332)(1210,3331)
	(1177,3331)(1141,3329)(1104,3328)
	(1066,3326)(1028,3324)(991,3321)
	(955,3318)(921,3314)(888,3309)
	(857,3303)(829,3297)(802,3289)
	(776,3280)(752,3270)(728,3258)
	(706,3246)(685,3232)(664,3218)
	(642,3202)(621,3185)(599,3167)
	(577,3148)(555,3128)(533,3107)
	(511,3086)(489,3064)(466,3042)
	(444,3020)(423,2998)(401,2975)
	(380,2952)(359,2930)(338,2907)
	(318,2885)(299,2863)(281,2841)
	(263,2819)(247,2797)(231,2776)
	(217,2755)(203,2733)(190,2709)
	(177,2685)(166,2659)(156,2632)
	(146,2604)(137,2573)(128,2540)
	(119,2506)(111,2470)(104,2433)
	(96,2395)(89,2357)(83,2320)
	(77,2284)(71,2251)(67,2222)
	(63,2196)(59,2175)(57,2159)
	(55,2147)(54,2139)(53,2135)(53,2133)
\path(353,2058)(353,2061)(352,2067)
	(351,2078)(349,2094)(347,2114)
	(344,2138)(342,2164)(340,2192)
	(338,2219)(338,2245)(338,2270)
	(339,2294)(342,2316)(347,2337)
	(353,2358)(359,2375)(366,2393)
	(375,2411)(384,2429)(394,2448)
	(405,2468)(417,2488)(429,2508)
	(442,2529)(455,2549)(468,2569)
	(481,2590)(494,2610)(507,2629)
	(519,2648)(532,2667)(544,2684)
	(555,2701)(567,2718)(578,2733)
	(592,2751)(606,2769)(620,2785)
	(633,2801)(647,2817)(660,2833)
	(674,2848)(687,2863)(700,2878)
	(714,2892)(728,2905)(742,2918)
	(756,2930)(771,2941)(787,2950)
	(803,2958)(820,2964)(838,2969)
	(858,2972)(881,2973)(905,2973)
	(932,2973)(960,2971)(989,2969)
	(1017,2967)(1043,2964)(1065,2962)
	(1081,2960)(1093,2959)(1100,2958)(1103,2958)
\path(1553,3183)(1555,3183)(1558,3183)
	(1564,3183)(1575,3184)(1589,3184)
	(1607,3184)(1631,3185)(1658,3186)
	(1690,3186)(1726,3187)(1766,3188)
	(1808,3189)(1852,3190)(1898,3190)
	(1944,3191)(1991,3192)(2037,3193)
	(2082,3193)(2125,3193)(2167,3193)
	(2206,3193)(2244,3193)(2279,3192)
	(2313,3191)(2344,3190)(2373,3189)
	(2401,3187)(2428,3185)(2453,3183)
	(2484,3180)(2514,3176)(2543,3172)
	(2571,3167)(2599,3162)(2627,3157)
	(2654,3152)(2680,3146)(2707,3140)
	(2733,3134)(2759,3128)(2785,3122)
	(2810,3116)(2835,3109)(2860,3103)
	(2884,3096)(2908,3089)(2931,3082)
	(2953,3075)(2975,3067)(2996,3059)
	(3016,3051)(3035,3042)(3053,3033)
	(3075,3020)(3094,3006)(3112,2991)
	(3128,2976)(3143,2960)(3156,2943)
	(3168,2926)(3179,2908)(3190,2890)
	(3200,2872)(3211,2854)(3221,2836)
	(3231,2818)(3240,2801)(3250,2784)
	(3260,2767)(3269,2750)(3278,2733)
	(3287,2714)(3295,2694)(3302,2673)
	(3309,2650)(3315,2626)(3321,2599)
	(3327,2571)(3333,2543)(3338,2516)
	(3343,2491)(3347,2470)(3349,2454)
	(3351,2442)(3352,2436)(3353,2433)
\path(1253,2883)(1255,2883)(1259,2883)
	(1267,2883)(1280,2883)(1297,2883)
	(1320,2883)(1348,2883)(1380,2883)
	(1416,2883)(1456,2883)(1499,2883)
	(1543,2883)(1588,2883)(1633,2883)
	(1677,2883)(1719,2883)(1760,2883)
	(1798,2883)(1834,2883)(1868,2883)
	(1899,2883)(1928,2883)(1954,2883)
	(1979,2883)(2003,2883)(2032,2883)
	(2060,2883)(2087,2883)(2112,2883)
	(2136,2884)(2159,2884)(2182,2885)
	(2203,2885)(2225,2886)(2246,2887)
	(2266,2887)(2287,2888)(2307,2888)
	(2328,2888)(2349,2888)(2369,2888)
	(2390,2887)(2411,2886)(2432,2885)
	(2453,2883)(2476,2881)(2498,2878)
	(2520,2875)(2542,2872)(2564,2869)
	(2586,2866)(2608,2863)(2629,2859)
	(2651,2856)(2672,2852)(2694,2849)
	(2714,2844)(2735,2840)(2755,2835)
	(2774,2829)(2793,2823)(2811,2816)
	(2828,2808)(2846,2798)(2864,2787)
	(2880,2775)(2896,2762)(2912,2748)
	(2928,2733)(2943,2718)(2958,2703)
	(2973,2687)(2987,2671)(3000,2656)
	(3013,2641)(3025,2626)(3036,2611)
	(3045,2597)(3053,2583)(3060,2567)
	(3065,2550)(3067,2531)(3068,2511)
	(3068,2488)(3067,2465)(3064,2440)
	(3062,2417)(3059,2397)(3056,2380)
	(3055,2368)(3054,2361)(3053,2358)
\path(653,1533)(353,1233)
\end{picture}
}

%% file: draws/crossing.tex
\begingroup\makeatletter\ifx\SetFigFont\undefined
\def\x#1#2#3#4#5#6#7\relax{\def\x{#1#2#3#4#5#6}}%
\expandafter\x\fmtname xxxxxx\relax \def\y{splain}%
\ifx\x\y   
\gdef\SetFigFont#1#2#3{%
  \ifnum #1<17\tiny\else \ifnum #1<20\small\else
  \ifnum #1<24\normalsize\else \ifnum #1<29\large\else
  \ifnum #1<34\Large\else \ifnum #1<41\LARGE\else
     \huge\fi\fi\fi\fi\fi\fi
  \csname #3\endcsname}%
\else
\gdef\SetFigFont#1#2#3{\begingroup
  \count@#1\relax \ifnum 25<\count@\count@25\fi
  \def\x{\endgroup\@setsize\SetFigFont{#2pt}}%
  \expandafter\x
    \csname \romannumeral\the\count@ pt\expandafter\endcsname
    \csname @\romannumeral\the\count@ pt\endcsname
  \csname #3\endcsname}%
\fi
\fi\endgroup
\begin{picture}(5120,1560)(0,-10)
\thicklines
\path(4334.618,1516.728)(4212.000,1533.000)(4312.535,1460.940)
\put(4033.875,1083.000){\arc{967.943}{5.0893}{7.4771}}
\path(5011.465,1460.940)(5112.000,1533.000)(4989.382,1516.728)
\put(5290.125,1083.000){\arc{967.943}{1.9477}{4.3355}}
\path(12,633)(912,1533)
\path(848.360,1426.934)(912.000,1533.000)(805.934,1469.360)
\path(118.066,1469.360)(12.000,1533.000)(75.640,1426.934)
\path(12,1533)(387,1158)
\path(2218.066,1469.360)(2112.000,1533.000)(2175.640,1426.934)
\path(2112,1533)(3012,633)
\path(2948.360,1426.934)(3012.000,1533.000)(2905.934,1469.360)
\path(3012,1533)(2637,1158)
\path(537,1008)(912,633)
\path(2487,1008)(2112,633)
\put(237,33){\makebox(0,0)[lb]{$L^+$}}
\put(4437,33){\makebox(0,0)[lb]{$L^0$}}
\put(2337,33){\makebox(0,0)[lb]{$L^-$}}
\end{picture}

%% file: draws/Dehn.tex
\begingroup\makeatletter\ifx\SetFigFont\undefined%
\gdef\SetFigFont#1#2#3#4#5{%
  \reset@font\fontsize{#1}{#2pt}%
  \fontfamily{#3}\fontseries{#4}\fontshape{#5}%
  \selectfont}%
\fi\endgroup%
{\renewcommand{\dashlinestretch}{30}
\begin{picture}(2950,927)(0,-10)
\path(238,900)(463,675)
\path(238,300)(838,900)
\path(774.360,793.934)(838.000,900.000)(731.934,836.360)
\path(613,525)(838,300)
\path(731.934,363.640)(838.000,300.000)(774.360,406.066)
\path(238,600)(235,598)(229,593)
	(218,584)(203,572)(184,556)
	(161,538)(137,517)(112,496)
	(89,476)(68,456)(50,437)
	(35,420)(24,404)(17,389)
	(13,375)(12,362)(13,350)
	(15,337)(18,323)(23,310)
	(28,295)(33,281)(39,266)
	(46,252)(53,238)(61,224)
	(71,210)(81,197)(93,186)
	(107,175)(123,165)(142,157)
	(163,150)(181,146)(201,142)
	(223,139)(248,137)(274,135)
	(302,133)(331,132)(363,131)
	(395,130)(429,129)(463,129)
	(497,129)(531,129)(566,129)
	(599,130)(632,131)(663,132)
	(694,133)(722,135)(749,137)
	(774,139)(797,142)(818,146)
	(838,150)(864,158)(888,168)
	(909,179)(928,192)(945,207)
	(961,222)(975,239)(989,255)
	(1002,272)(1014,289)(1026,306)
	(1036,321)(1045,336)(1053,350)
	(1059,363)(1063,375)(1064,388)
	(1060,401)(1052,414)(1040,428)
	(1025,443)(1007,458)(987,474)
	(967,488)(949,501)(933,511)
	(923,519)(916,523)(913,525)
\path(2563,525)(2338,300)
\path(2338,900)(2938,300)
\path(2831.934,363.640)(2938.000,300.000)(2874.360,406.066)
\path(2874.360,793.934)(2938.000,900.000)(2831.934,836.360)
\path(2938,900)(2713,675)
\put(1213,0){\makebox(0,0)[lb]{$+1$}}
\put(1438,525){\makebox(0,0)[lb]{$=$}}
\end{picture}
}

%% file: draws/step1.tex
\begingroup\makeatletter\ifx\SetFigFont\undefined%
\gdef\SetFigFont#1#2#3#4#5{%
  \reset@font\fontsize{#1}{#2pt}%
  \fontfamily{#3}\fontseries{#4}\fontshape{#5}%
  \selectfont}%
\fi\endgroup%
{\renewcommand{\dashlinestretch}{30}
\begin{picture}(6912,1849)(0,-10)
\put(5992.819,935.553){\arc{1214.869}{3.3614}{8.8747}}
\put(5999.472,911.662){\arc{1806.284}{3.2318}{9.0043}}
\path(4950,843)(4950,693)
\path(5250,843)(5250,693)
\path(5475,1218)(5478,1216)(5483,1213)
	(5493,1208)(5507,1199)(5525,1189)
	(5546,1177)(5569,1163)(5592,1149)
	(5614,1135)(5635,1122)(5653,1109)
	(5669,1098)(5682,1087)(5692,1077)
	(5700,1068)(5707,1057)(5712,1046)
	(5716,1035)(5718,1024)(5720,1013)
	(5721,1001)(5721,989)(5721,978)
	(5720,967)(5718,956)(5716,945)
	(5712,936)(5707,926)(5700,918)
	(5691,910)(5678,903)(5662,896)
	(5643,889)(5620,881)(5595,874)
	(5569,867)(5543,860)(5520,854)
	(5501,849)(5487,846)(5479,844)
	(5476,843)(5475,843)
\path(5550,1368)(5552,1366)(5558,1363)
	(5567,1356)(5581,1346)(5599,1334)
	(5620,1318)(5645,1301)(5671,1282)
	(5698,1262)(5724,1243)(5749,1225)
	(5772,1208)(5792,1192)(5810,1178)
	(5826,1165)(5839,1154)(5850,1143)
	(5862,1130)(5873,1118)(5882,1107)
	(5890,1096)(5897,1086)(5903,1077)
	(5909,1068)(5914,1059)(5918,1050)
	(5922,1040)(5925,1029)(5926,1018)
	(5927,1006)(5925,993)(5922,981)
	(5917,968)(5911,955)(5905,940)
	(5897,925)(5889,909)(5880,893)
	(5871,877)(5861,861)(5851,845)
	(5841,830)(5830,815)(5818,802)
	(5805,789)(5791,778)(5775,768)
	(5760,760)(5742,753)(5722,746)
	(5699,740)(5674,734)(5645,729)
	(5614,723)(5582,718)(5549,713)
	(5516,708)(5486,704)(5459,701)
	(5437,698)(5421,696)(5409,694)
	(5403,693)(5400,693)
\path(5400,693)(5398,693)(5392,693)
	(5383,693)(5368,692)(5348,692)
	(5322,692)(5290,691)(5254,691)
	(5214,690)(5171,690)(5126,689)
	(5081,689)(5036,688)(4992,688)
	(4951,688)(4912,688)(4876,688)
	(4843,688)(4813,689)(4787,689)
	(4764,690)(4743,692)(4725,693)
	(4698,696)(4676,700)(4659,704)
	(4645,709)(4634,715)(4625,721)
	(4618,727)(4612,733)(4606,739)
	(4600,746)(4594,751)(4587,757)
	(4581,763)(4575,768)(4569,776)
	(4565,783)(4564,792)(4564,800)
	(4564,809)(4566,818)(4568,826)
	(4570,833)(4572,839)(4575,843)
	(4581,847)(4591,848)(4604,848)
	(4619,847)(4633,845)(4644,844)
	(4649,843)(4650,843)
\path(5175,1293)(5172,1294)(5165,1296)
	(5153,1300)(5135,1305)(5113,1312)
	(5086,1320)(5058,1328)(5028,1336)
	(4999,1344)(4972,1351)(4948,1357)
	(4926,1362)(4906,1365)(4890,1367)
	(4875,1368)(4857,1367)(4841,1364)
	(4827,1360)(4814,1354)(4801,1348)
	(4789,1341)(4778,1334)(4767,1326)
	(4756,1318)(4745,1310)(4735,1302)
	(4725,1293)(4716,1283)(4708,1271)
	(4699,1257)(4691,1240)(4682,1221)
	(4674,1202)(4666,1183)(4659,1166)
	(4654,1154)(4651,1147)(4650,1144)(4650,1143)
\path(4650,1143)(4651,1139)(4654,1132)
	(4659,1119)(4665,1100)(4673,1078)
	(4681,1053)(4690,1027)(4698,1003)
	(4706,980)(4712,960)(4718,944)
	(4722,930)(4725,918)(4729,902)
	(4731,889)(4732,880)(4732,873)
	(4732,867)(4731,860)(4729,852)
	(4725,843)(4720,835)(4714,828)
	(4706,819)(4697,811)(4688,802)
	(4678,793)(4669,785)(4661,778)
	(4655,772)(4650,768)(4649,767)
	(4648,766)(4647,765)(4646,764)
	(4645,763)(4646,764)(4647,765)
	(4648,766)(4649,767)(4650,768)
\path(4500,843)(4500,846)(4500,853)
	(4500,863)(4500,877)(4500,894)
	(4500,912)(4500,930)(4500,947)
	(4500,963)(4500,978)(4500,993)
	(4500,1007)(4500,1022)(4500,1040)
	(4500,1060)(4500,1083)(4500,1107)
	(4500,1132)(4500,1156)(4500,1177)
	(4500,1195)(4500,1207)(4500,1215)(4500,1218)
\path(4500,1218)(4501,1221)(4504,1228)
	(4509,1240)(4516,1255)(4524,1273)
	(4532,1292)(4541,1310)(4549,1327)
	(4558,1342)(4566,1355)(4575,1368)
	(4583,1379)(4592,1390)(4602,1401)
	(4612,1413)(4622,1425)(4632,1438)
	(4643,1450)(4654,1462)(4665,1474)
	(4676,1485)(4688,1495)(4700,1504)
	(4712,1512)(4725,1518)(4739,1523)
	(4754,1526)(4770,1527)(4786,1528)
	(4804,1528)(4822,1527)(4841,1525)
	(4859,1523)(4877,1522)(4895,1520)
	(4911,1519)(4925,1518)(4938,1518)
	(4950,1518)(4966,1518)(4978,1518)
	(4986,1518)(4991,1518)(4996,1518)
	(5003,1518)(5012,1518)(5025,1518)
	(5035,1518)(5046,1518)(5058,1519)
	(5072,1520)(5086,1521)(5100,1522)
	(5114,1522)(5128,1523)(5142,1523)
	(5154,1522)(5165,1521)(5175,1518)
	(5186,1513)(5196,1507)(5206,1498)
	(5216,1486)(5226,1474)(5236,1462)
	(5243,1453)(5248,1446)(5250,1443)
\path(4800,843)(4803,843)(4809,843)
	(4819,843)(4833,843)(4851,843)
	(4871,843)(4893,843)(4915,843)
	(4937,843)(4959,843)(4980,843)
	(5002,843)(5025,843)(5043,843)
	(5062,843)(5083,843)(5107,843)
	(5134,843)(5163,843)(5195,843)
	(5228,843)(5263,843)(5299,843)
	(5333,843)(5366,843)(5396,843)
	(5421,843)(5441,843)(5457,843)
	(5467,843)(5472,843)(5475,843)
\put(1715.625,918.000){\arc{1218.750}{3.3903}{9.1761}}
\put(1718.750,918.000){\arc{1812.500}{3.3079}{9.2585}}
\put(975.000,1405.500){\arc{1275.000}{1.0808}{2.0608}}
\path(675,1293)(825,1293)
\path(1275,1293)(1276,1293)(1279,1290)
	(1287,1283)(1301,1272)(1320,1256)
	(1342,1237)(1365,1216)(1389,1195)
	(1410,1175)(1428,1156)(1443,1138)
	(1455,1122)(1464,1108)(1470,1094)
	(1474,1081)(1475,1068)(1474,1055)
	(1470,1042)(1464,1028)(1455,1014)
	(1443,998)(1428,980)(1410,961)
	(1389,941)(1365,920)(1342,899)
	(1320,880)(1301,864)(1287,853)
	(1279,846)(1276,843)(1275,843)
\path(675,1293)(674,1293)(671,1290)
	(663,1283)(649,1272)(630,1256)
	(608,1237)(585,1216)(561,1195)
	(540,1175)(522,1156)(507,1138)
	(495,1122)(486,1108)(480,1094)
	(476,1081)(475,1068)(476,1055)
	(480,1042)(486,1028)(495,1014)
	(507,998)(522,980)(540,961)
	(561,941)(585,920)(608,899)
	(630,880)(649,864)(663,853)
	(671,846)(674,843)(675,843)
\path(4950,843)(5250,693)
\path(4950,693)(5250,843)
\path(2325,918)(2625,918)
\path(2250,1143)(2550,1293)
\path(2100,1368)(2325,1593)
\path(2250,693)(2550,543)
\path(2100,468)(2325,243)
\path(6600,918)(6900,918)
\path(6600,1143)(6825,1218)
\path(6525,1293)(6675,1443)
\path(6600,693)(6825,618)
\path(6450,543)(6675,318)
\put(0,843){\makebox(0,0)[lb]{$K$}}
\put(3600,843){\makebox(0,0)[lb]{$D_mK$}}
\end{picture}
}

%% file: draws/step2.tex
\begingroup\makeatletter\ifx\SetFigFont\undefined%
\gdef\SetFigFont#1#2#3#4#5{%
  \reset@font\fontsize{#1}{#2pt}%
  \fontfamily{#3}\fontseries{#4}\fontshape{#5}%
  \selectfont}%
\fi\endgroup%
{\renewcommand{\dashlinestretch}{30}
\begin{picture}(7962,1998)(0,-10)
\put(2017.819,935.553){\arc{1214.869}{3.3614}{8.8747}}
\put(2024.472,911.662){\arc{1806.284}{3.2318}{9.0043}}
\put(7049.473,986.662){\arc{1806.286}{3.2318}{9.0043}}
\put(7042.818,1010.553){\arc{1214.867}{3.3614}{8.8747}}
\put(5287.500,1105.500){\arc{237.171}{1.8925}{4.3906}}
\put(5737.500,1105.500){\arc{237.171}{5.0341}{7.5322}}
\path(975,843)(975,693)
\path(1275,843)(1275,693)
\path(6750,1143)(6975,1143)(6975,918)
	(6750,918)(6750,1143)
\path(975,843)(1275,693)
\path(975,693)(1275,843)
\path(6750,1143)(6975,918)
\path(6975,1143)(6750,918)
\path(2625,918)(2925,918)
\path(2625,1143)(2850,1218)
\path(2475,1368)(2700,1518)
\path(2625,693)(2850,618)
\path(2475,543)(2700,318)
\path(7650,1068)(7950,1068)
\path(7575,1293)(7800,1443)
\path(7425,1518)(7575,1668)
\path(7650,843)(7875,768)
\path(7575,693)(7800,468)
\path(1500,1218)(1503,1216)(1508,1213)
	(1518,1208)(1532,1199)(1550,1189)
	(1571,1177)(1594,1163)(1617,1149)
	(1639,1135)(1660,1122)(1678,1109)
	(1694,1098)(1707,1087)(1717,1077)
	(1725,1068)(1732,1057)(1737,1046)
	(1741,1035)(1743,1024)(1745,1013)
	(1746,1001)(1746,989)(1746,978)
	(1745,967)(1743,956)(1741,945)
	(1737,936)(1732,926)(1725,918)
	(1716,910)(1703,903)(1687,896)
	(1668,889)(1645,881)(1620,874)
	(1594,867)(1568,860)(1545,854)
	(1526,849)(1512,846)(1504,844)
	(1501,843)(1500,843)
\path(1575,1368)(1577,1366)(1583,1363)
	(1592,1356)(1606,1346)(1624,1334)
	(1645,1318)(1670,1301)(1696,1282)
	(1723,1262)(1749,1243)(1774,1225)
	(1797,1208)(1817,1192)(1835,1178)
	(1851,1165)(1864,1154)(1875,1143)
	(1887,1130)(1898,1118)(1907,1107)
	(1915,1096)(1922,1086)(1928,1077)
	(1934,1068)(1939,1059)(1943,1050)
	(1947,1040)(1950,1029)(1951,1018)
	(1952,1006)(1950,993)(1947,981)
	(1942,968)(1936,955)(1930,940)
	(1922,925)(1914,909)(1905,893)
	(1896,877)(1886,861)(1876,845)
	(1866,830)(1855,815)(1843,802)
	(1830,789)(1816,778)(1800,768)
	(1785,760)(1767,753)(1747,746)
	(1724,740)(1699,734)(1670,729)
	(1639,723)(1607,718)(1574,713)
	(1541,708)(1511,704)(1484,701)
	(1462,698)(1446,696)(1434,694)
	(1428,693)(1425,693)
\path(1200,1293)(1197,1294)(1190,1296)
	(1178,1300)(1160,1305)(1138,1312)
	(1111,1320)(1083,1328)(1053,1336)
	(1024,1344)(997,1351)(973,1357)
	(951,1362)(931,1365)(915,1367)
	(900,1368)(882,1367)(866,1364)
	(852,1360)(839,1354)(826,1348)
	(814,1341)(803,1334)(792,1326)
	(781,1318)(770,1310)(760,1302)
	(750,1293)(741,1283)(733,1271)
	(724,1257)(716,1240)(707,1221)
	(699,1202)(691,1183)(684,1166)
	(679,1154)(676,1147)(675,1144)(675,1143)
\path(525,1218)(526,1221)(529,1228)
	(534,1240)(541,1255)(549,1273)
	(557,1292)(566,1310)(574,1327)
	(583,1342)(591,1355)(600,1368)
	(608,1379)(617,1390)(627,1401)
	(637,1413)(647,1425)(657,1438)
	(668,1450)(679,1462)(690,1474)
	(701,1485)(713,1495)(725,1504)
	(737,1512)(750,1518)(764,1523)
	(779,1526)(795,1527)(811,1528)
	(829,1528)(847,1527)(866,1525)
	(884,1523)(902,1522)(920,1520)
	(936,1519)(950,1518)(963,1518)
	(975,1518)(991,1518)(1003,1518)
	(1011,1518)(1016,1518)(1021,1518)
	(1028,1518)(1037,1518)(1050,1518)
	(1060,1518)(1071,1518)(1083,1519)
	(1097,1520)(1111,1521)(1125,1522)
	(1139,1522)(1153,1523)(1167,1523)
	(1179,1522)(1190,1521)(1200,1518)
	(1211,1513)(1221,1507)(1231,1498)
	(1241,1486)(1251,1474)(1261,1462)
	(1268,1453)(1273,1446)(1275,1443)
\path(825,843)(828,843)(834,843)
	(844,843)(858,843)(876,843)
	(896,843)(918,843)(940,843)
	(962,843)(984,843)(1005,843)
	(1027,843)(1050,843)(1068,843)
	(1087,843)(1108,843)(1132,843)
	(1159,843)(1188,843)(1220,843)
	(1253,843)(1288,843)(1324,843)
	(1358,843)(1391,843)(1421,843)
	(1446,843)(1466,843)(1482,843)
	(1492,843)(1497,843)(1500,843)
\path(6525,1293)(6526,1293)(6529,1292)
	(6537,1290)(6551,1286)(6570,1280)
	(6593,1273)(6619,1265)(6645,1257)
	(6670,1250)(6693,1242)(6712,1235)
	(6728,1229)(6741,1223)(6750,1218)
	(6759,1210)(6764,1202)(6765,1192)
	(6764,1182)(6761,1170)(6757,1160)
	(6754,1151)(6751,1146)(6750,1143)
\path(6600,1443)(6602,1442)(6607,1440)
	(6616,1437)(6630,1433)(6648,1426)
	(6671,1418)(6697,1409)(6726,1399)
	(6756,1388)(6787,1377)(6817,1366)
	(6846,1355)(6873,1344)(6897,1334)
	(6918,1325)(6936,1316)(6952,1308)
	(6965,1301)(6975,1293)(6986,1282)
	(6994,1271)(6999,1258)(7001,1245)
	(7000,1230)(6998,1214)(6994,1198)
	(6989,1183)(6985,1169)(6981,1158)
	(6978,1150)(6976,1145)(6975,1143)
\path(6975,918)(6975,917)(6978,914)
	(6983,907)(6991,895)(7002,878)
	(7014,859)(7026,840)(7036,821)
	(7044,804)(7050,790)(7051,778)
	(7050,768)(7045,759)(7038,751)
	(7030,744)(7019,737)(7008,730)
	(6996,723)(6983,717)(6969,712)
	(6955,706)(6938,701)(6920,697)
	(6900,693)(6885,691)(6867,690)
	(6847,689)(6824,688)(6799,688)
	(6770,688)(6739,688)(6707,689)
	(6674,689)(6641,690)(6611,691)
	(6584,691)(6562,692)(6546,692)
	(6534,693)(6528,693)(6525,693)
\path(6750,918)(6750,917)(6751,916)
	(6752,914)(6753,911)(6754,908)
	(6755,905)(6756,900)(6756,895)
	(6756,890)(6756,885)(6754,880)
	(6752,875)(6749,870)(6744,866)
	(6738,862)(6729,857)(6719,854)
	(6707,850)(6692,846)(6675,843)
	(6658,840)(6639,837)(6618,834)
	(6595,831)(6569,828)(6542,824)
	(6513,821)(6482,817)(6449,814)
	(6416,810)(6381,806)(6346,802)
	(6311,798)(6275,795)(6240,791)
	(6205,788)(6170,784)(6137,781)
	(6105,778)(6075,776)(6046,773)
	(6018,771)(5992,770)(5968,769)
	(5946,768)(5925,768)(5894,769)
	(5866,770)(5840,773)(5817,777)
	(5796,781)(5776,787)(5758,792)
	(5740,798)(5723,805)(5707,811)
	(5691,817)(5676,823)(5661,828)
	(5648,833)(5635,838)(5625,843)
	(5616,848)(5610,854)(5606,860)
	(5604,867)(5605,875)(5607,882)
	(5610,890)(5614,898)(5617,905)
	(5620,911)(5623,915)(5624,917)(5625,918)
\path(6525,693)(6522,692)(6516,691)
	(6506,689)(6490,685)(6469,680)
	(6443,675)(6415,668)(6384,662)
	(6352,655)(6320,648)(6290,642)
	(6262,636)(6235,631)(6211,627)
	(6189,623)(6169,620)(6150,618)
	(6129,616)(6109,614)(6089,614)
	(6069,613)(6050,613)(6031,613)
	(6012,614)(5993,614)(5974,615)
	(5956,616)(5937,617)(5919,617)
	(5901,618)(5884,618)(5867,618)
	(5850,618)(5832,618)(5814,618)
	(5797,617)(5781,616)(5765,614)
	(5749,613)(5734,611)(5719,609)
	(5704,608)(5688,608)(5673,609)
	(5657,610)(5641,613)(5625,618)
	(5611,623)(5597,630)(5582,638)
	(5566,647)(5551,658)(5535,669)
	(5518,681)(5502,693)(5486,705)
	(5470,717)(5455,728)(5441,739)
	(5429,748)(5418,756)(5408,763)
	(5400,768)(5397,770)(5395,771)
	(5393,773)(5391,774)(5389,775)
	(5388,776)(5387,776)(5386,776)
	(5385,776)(5385,775)(5385,774)
	(5386,773)(5386,772)(5387,771)
	(5388,770)(5389,769)(5389,768)
	(5390,767)(5391,766)(5392,765)
	(5393,764)(5394,763)(5395,762)
	(5396,761)(5396,760)(5397,760)
	(5398,760)(5399,760)(5400,761)
	(5400,762)(5400,763)(5400,765)
	(5400,766)(5400,768)(5400,775)
	(5400,786)(5400,799)(5400,816)
	(5400,835)(5400,855)(5400,875)
	(5400,893)(5400,906)(5400,914)
	(5400,917)(5400,918)
\path(5250,993)(5253,992)(5259,988)
	(5270,983)(5285,976)(5304,968)
	(5325,958)(5348,948)(5370,940)
	(5392,932)(5413,926)(5434,921)
	(5454,919)(5475,918)(5492,918)
	(5511,920)(5532,923)(5556,926)
	(5582,931)(5610,937)(5641,943)
	(5673,951)(5705,958)(5737,965)
	(5766,972)(5792,979)(5814,984)
	(5830,988)(5841,991)(5847,992)(5850,993)
\path(6225,1368)(6222,1369)(6216,1371)
	(6204,1375)(6188,1380)(6167,1387)
	(6142,1395)(6115,1403)(6087,1411)
	(6059,1419)(6032,1426)(6008,1432)
	(5985,1437)(5964,1440)(5944,1442)
	(5925,1443)(5906,1442)(5887,1441)
	(5867,1439)(5847,1436)(5826,1432)
	(5804,1428)(5783,1424)(5761,1420)
	(5739,1415)(5719,1410)(5699,1404)
	(5680,1398)(5663,1392)(5648,1385)
	(5636,1377)(5625,1368)(5616,1356)
	(5610,1342)(5606,1326)(5604,1306)
	(5605,1283)(5607,1259)(5610,1233)
	(5614,1208)(5617,1186)(5620,1167)
	(5623,1154)(5624,1147)(5625,1143)
\path(5400,1143)(5397,1146)(5392,1153)
	(5384,1163)(5373,1177)(5361,1194)
	(5349,1212)(5339,1230)(5331,1247)
	(5325,1263)(5324,1278)(5325,1293)
	(5328,1305)(5333,1318)(5339,1331)
	(5346,1345)(5354,1360)(5363,1375)
	(5373,1390)(5382,1405)(5393,1421)
	(5403,1436)(5414,1451)(5425,1466)
	(5437,1480)(5449,1493)(5461,1506)
	(5475,1518)(5488,1528)(5502,1538)
	(5517,1548)(5533,1557)(5549,1567)
	(5567,1576)(5585,1585)(5603,1594)
	(5621,1604)(5640,1612)(5659,1621)
	(5677,1629)(5695,1637)(5712,1645)
	(5729,1652)(5745,1658)(5760,1663)
	(5775,1668)(5793,1673)(5811,1676)
	(5828,1678)(5844,1679)(5860,1679)
	(5876,1679)(5891,1679)(5906,1678)
	(5921,1676)(5937,1675)(5952,1674)
	(5968,1672)(5984,1670)(6000,1668)
	(6016,1665)(6032,1661)(6049,1656)
	(6066,1651)(6084,1644)(6102,1637)
	(6120,1630)(6137,1624)(6154,1617)
	(6171,1610)(6186,1605)(6201,1600)
	(6214,1596)(6225,1593)(6242,1589)
	(6256,1588)(6268,1588)(6280,1589)
	(6290,1591)(6296,1592)(6299,1593)(6300,1593)
\path(600,1068)(599,1068)(596,1069)
	(587,1070)(573,1073)(555,1076)
	(534,1079)(512,1082)(493,1083)
	(476,1083)(463,1081)(455,1076)
	(450,1068)(448,1059)(449,1047)
	(451,1033)(455,1017)(460,999)
	(467,979)(473,959)(481,939)
	(488,919)(495,900)(503,883)
	(510,868)(517,855)(525,843)
	(535,831)(546,822)(558,814)
	(572,807)(586,800)(600,795)
	(614,790)(628,785)(642,780)
	(654,776)(665,772)(675,768)
	(688,764)(701,763)(713,763)
	(726,764)(738,766)(746,767)
	(749,768)(750,768)
\path(825,843)(675,843)
\path(600,693)(603,691)(610,686)
	(620,679)(634,669)(651,658)
	(669,647)(687,637)(704,629)
	(720,623)(735,619)(750,618)
	(764,619)(778,621)(792,625)
	(806,630)(820,637)(834,644)
	(848,652)(863,660)(878,667)
	(895,674)(912,681)(931,686)
	(952,690)(975,693)(993,695)
	(1012,696)(1033,697)(1057,698)
	(1084,698)(1113,698)(1145,698)
	(1178,698)(1213,697)(1249,697)
	(1283,696)(1316,695)(1346,695)
	(1371,694)(1391,694)(1407,693)
	(1417,693)(1422,693)(1425,693)
\path(675,1143)(675,1142)(675,1139)
	(676,1131)(677,1116)(679,1096)
	(681,1072)(683,1045)(684,1018)
	(685,993)(685,970)(685,951)
	(683,936)(680,925)(675,918)
	(668,914)(657,914)(644,919)
	(627,927)(608,938)(588,950)
	(568,964)(550,975)(537,985)
	(529,990)(526,993)(525,993)
\put(150,1668){\makebox(0,0)[lb]{$D_mK$}}
\put(5625,1818){\makebox(0,0)[lb]{$C$}}
\put(0,768){\makebox(0,0)[lb]{$C$}}
\put(4350,993){\makebox(0,0)[lb]{$D_mK$}}
\end{picture}
}

%% file: draws/step4.tex
\begingroup\makeatletter\ifx\SetFigFont\undefined%
\gdef\SetFigFont#1#2#3#4#5{%
  \reset@font\fontsize{#1}{#2pt}%
  \fontfamily{#3}\fontseries{#4}\fontshape{#5}%
  \selectfont}%
\fi\endgroup%
{\renewcommand{\dashlinestretch}{30}
\begin{picture}(7528,7242)(0,-10)
\put(1965.000,5263.000){\arc{2121.320}{2.9997}{4.8543}}
\put(1065.000,3605.500){\arc{1515.000}{1.7701}{7.6546}}
\put(1065.000,3607.643){\arc{2110.715}{1.7134}{7.7114}}
\put(1965.000,1963.000){\arc{2121.320}{1.4289}{3.2835}}
\put(1991.786,5229.072){\arc{1570.819}{2.9933}{4.7738}}
\put(1991.786,1996.929){\arc{1570.820}{1.5094}{3.2899}}
\put(3771.000,1063.000){\arc{2109.442}{0.1427}{6.1405}}
\put(3772.500,1063.000){\arc{1515.000}{0.1993}{6.0838}}
\put(5539.000,1997.000){\arc{1571.060}{6.1348}{7.9164}}
\put(5565.000,1963.000){\arc{2121.320}{6.1413}{7.9959}}
\put(5565.000,5263.000){\arc{2121.320}{4.5705}{6.4251}}
\put(6465.000,3470.000){\arc{1515.980}{4.9116}{10.7964}}
\put(6465.000,3469.000){\arc{2109.442}{4.8551}{10.8529}}
\put(3757.000,6163.000){\arc{1514.020}{3.3411}{9.2253}}
\put(3759.000,6163.000){\arc{2109.442}{3.2843}{9.2821}}
\put(5532.000,5236.000){\arc{1568.940}{4.5629}{6.3444}}
\path(1215,5113)(1215,4663)
\path(915,5113)(915,4663)
\path(915,4363)(915,2113)
\path(1215,4363)(1215,2113)
\path(2040,1213)(2715,1213)
\path(2115,913)(2715,913)
\path(3015,1213)(5490,1213)
\path(3015,913)(5415,913)
\path(2115,6313)(4515,6313)
\path(2040,6013)(4515,6013)
\path(4815,6313)(5415,6313)
\path(4815,6013)(5415,6013)
\path(6315,5188)(6315,2713)
\path(6615,5113)(6615,2713)
\path(6315,2413)(6315,2038)
\path(6615,2413)(6615,2038)
\path(915,3913)(1215,3913)
\path(915,3313)(1215,3313)
\path(3615,1213)(3615,913)
\path(4215,1213)(4215,913)
\path(6315,3913)(6615,3913)
\path(6315,3313)(6615,3313)
\path(3315,6313)(3315,6013)
\path(3915,6313)(3915,6013)
\path(240,4288)(465,4063)
\path(240,2938)(465,3163)
\path(3315,463)(3090,238)
\path(4215,388)(4365,238)
\path(7215,4213)(6990,3988)
\path(6915,2863)(7140,2638)
\path(4365,7063)(4140,6763)
\path(3015,6913)(3240,6688)
\path(915,3913)(1215,3313)
\path(1215,3913)(915,3313)
\path(3615,1213)(4215,913)
\path(3615,913)(4215,1213)
\path(6315,3913)(6615,3313)
\path(6615,3913)(6315,3313)
\path(3315,6313)(3915,6013)
\path(3315,6013)(3915,6313)
\path(90,3988)(315,3838)
\path(15,3613)(315,3613)
\path(90,3238)(315,3388)
\path(3765,313)(3765,13)
\path(3990,313)(4140,13)
\path(3540,313)(3390,88)
\path(7215,3463)(7515,3463)
\path(7140,3763)(7440,3838)
\path(7140,3163)(7440,2938)
\path(3615,7213)(3615,6913)
\path(3990,7213)(3915,6913)
\path(3240,7063)(3390,6838)
\end{picture}
}